\documentclass[centertags]{article}

\usepackage{amsmath}
\usepackage{amsthm}
\usepackage{amscd}
\usepackage{amssymb}
\usepackage{verbatim}
\usepackage[usenames]{color}

\title{Groups of homeomorphisms of one-manifolds, III: Nilpotent 
subgroups}
\author{Benson Farb\thanks{Supported in 
part by the NSF and by the Sloan Foundation.}\ \  and 
John Franks\thanks{Supported in part by NSF grants DMS9803346 and DMS0099640.}}

\newtheorem{theorem}{Theorem}[section]

\newtheorem{remark}[theorem]{Remark}

\newtheorem{thm}{Theorem}[section] 

\newtheorem{lem}[thm]{Lemma} 
\newtheorem{prop}[thm]{Proposition}
\newtheorem{defn}[thm]{Definition}

\newcommand{\R}{{\bf R}}
\newcommand{\Z}{\bf Z}
\newcommand{\T}{S}
\newcommand{\Q}{\bf Q}
\newcommand{\I}{I}
\newcommand{\N}{\mathcal N}
\def\endproof{$\diamond$ \bigskip}
\DeclareMathOperator{\Fix}{Fix}
\DeclareMathOperator{\Int}{Int}

\DeclareMathOperator{\Homeo}{Homeo}

\DeclareMathOperator{\Diff}{Diff}
\DeclareMathOperator{\PL}{PL}
\begin{document}
\maketitle
\begin{abstract}
This self-contained paper is part of a series \cite{FF1,FF2} seeking to
understand groups of homeomorphisms of manifolds in analogy with the
theory of Lie groups and their discrete subgroups.  Plante-Thurston
proved that every nilpotent subgroup of $\Diff^2(S^1)$ is abelian.  One
of our main results is a sharp converse: $\Diff^1(S^1)$
contains every finitely-generated, torsion-free nilpotent group.
\end{abstract}

\section{Introduction}

A basic aspect of the theory of linear groups 
is the structure of nilpotent groups.  In this paper we consider
nilpotent subgroups of $\Homeo(M)$ and $\Diff^r(M)$, where $M$ is 
the line $\R$, the circle $\T^1$, or the interval $I=[0,1]$. 
As we will see, the structure theory 
depends dramatically on the degree $r$ of regularity 
as well as on the topology of $M$.  

Throughout this paper all homeomorphisms will be orientation-preserving,
and all groups will consist only of such homeomorphisms.
Plante-Thurston \cite{PT} discovered that $C^2$ regularity imposes 
a severe restriction on nilpotent groups of diffeomorphisms.

\begin{theorem}
\label{theorem:C2interval}
\label{theorem:C^2_T^1}
Any nilpotent subgroup of $\Diff^2(\I),\ \Diff^2([0,1))$ or
$\Diff^2(\T^1)$ must be abelian.
\end{theorem}

\noindent
{\bf Remark. } In the case of subgroups of $\Diff^2(\I)$ and
$\Diff^2([0,1))$ this result was first proved by Plante and Thurston in
\cite{PT}.  They also proved that the group is virtually abelian in
the $\T^1$ case.  For completeness of exposition (and because it is
simple) we present a proof of the Plante-Thurston result about $\I$.
These results are related to a result of E. Ghys \cite{Gh}, who proved
that any solvable subgroup of $\Diff^\omega(\T^1)$ is metabelian.  In
\S\ref{section:PL} we prove the PL version of Theorem
\ref{theorem:C2interval}.
\bigskip

Our main result is that lowering the regularity from $C^2$ to $C^1$
produces a sharply contrasting situation, where every possibility can
occur. 

\begin{theorem}
\label{thm:fg.nilpotent}
Let $M=\R,\T^1,$ or $\I$.  Then every 
finitely-generated, torsion-free nilpotent group is isomorphic
to a subgroup of $\Diff^1(M)$.
\end{theorem}

Witte has observed (\cite{Wi}, Lemma 2.2) that the groups of
orientation-preserving homeomorphisms of $\R$ are precisely the
right-orderable groups.  Since torsion-free nilpotent groups are
right-orderable (see, e.g.\ \cite{MR}, p.37), it follows that such
groups are subgroups of $\Homeo(\R)$.  The content of Theorem 
\ref{thm:fg.nilpotent} is the increase in regularity from $C^0$ to the 
sharp regularity $C^1$.

\medskip
The situation for $\R$ is more complicated.  On the one hand there is no 
limit to the degree of nilpotence, even when regularity is high:

\begin{theorem}
\label{theorem:smooth:example}
$\Diff^{\infty}(\R)$ contains nilpotent subgroups of every degree of 
nilpotency.
\end{theorem}

On the other hand, even with just a little regularity, the {\em derived 
length} of nilpotent groups is greatly restricted.  

\begin{theorem}
\label{thm:derived.2}
Every nilpotent 
subgroup of $\Diff^2(\R)$ is metabelian, i.e.\ has abelian 
commutator subgroup.
\end{theorem}  

In particular, Theorem \ref{thm:derived.2} gives that the group 
of $n\times n$ upper-triangular with ones on the diagonal, while
admitting an effective action on $\R$ by $C^1$ diffeomorphisms, admits no
effective $C^2$ action on $\R$ if $n>3$.

Closely related to this algebraic restriction on nilpotent
groups which act smoothly is a topological restriction.

\begin{theorem}
\label{thm:fixed.abelian}
If $N$ is a nilpotent
subgroup of $\Diff^2(\R)$ and every element of $N$ has a fixed
point then $N$ is abelian.
\end{theorem}  

One open problem is to extend the above theory to solvable
subgroups, as well as to higher-dimensional manifolds.  The paper of
J. Plante \cite{P} contains a number of interesting results and
examples concerning solvable groups acting on $\R.$ Another
problem is to understand what happens between the degrees of
regularity $r=0$ and $r\geq 2$, where vastly differing phenomena
occur.

\bigskip
\noindent
{\bf Residually nilpotent groups. } 
A variation of our construction of actions of nilpotent groups can be 
used to construct actions of a much wider class of groups, the {\em
residually torsion-free nilpotent} groups, i.e.\ those groups where the
intersection of all terms in the lower central series is trivial.  
In \S\ref{section:residual} we prove the following result.

\begin{theorem}
\label{theorem:residual}
Let $M=\R,\T^1$ or $\I$.  Then $\Diff_+^1(M)$ contains every
finitely generated, residually torsion-free nilpotent group.
\end{theorem}

The class of finitely generated, residually torsion-free nilpotent
groups includes free groups, surface groups, and the Torelli groups, as
well as products of these groups.  

\bigskip
\noindent
{\bf Acknowledgements. } It is a pleasure to thank Dave Witte for making 
several useful suggestions, in particular pointing us to the references 
\cite{Ra}, \cite{P}, and \cite{MR}.

\section{Constructions}

In this section we construct examples of nilpotent groups acting on 
one-manifolds.

\subsection{An elementary observation}

There several obvious relations among actions on the three
spaces $\R,\T^1,$ and $\I$.  
Restricting an action on $\I$ to the interior of $\I$ gives an
action on $\R$.  One can also start with an action on $\R$ and
consider the action on the one-point compactification $\T^1$ or the
two-point compactification $\I$.  But this usually entails loss of
regularity of the action.  That is, a smooth or PL action on $\R$ will
generally give only an action by homeomorpisms on $\I$ or $T^1.$

\subsection{\boldmath{$C^1$} actions on $M$}

Let $\N_n$ denote the group of $n\times n$ lower-triangular integer
matrices with ones on the diagonal.  Our first main goal is to prove
that $\N_n$ admits $C^1$ actions on the real line. 

The proof is somewhat technical and so we will describe the strategy
before engaging in the details.  The traditional proof that $\N_n$
acts by homeomorphisms on the interval uses the fact that it is an
ordered group.  Suppose a group acts effectively on a countable
ordered set in a way that preserves the order (e.g. the set might be
the group itself), One can then produce an action by homeomorphisms on
$\I$ by embedding the countable set in an order preserving way in $\I$
and canonically extending the action of each group element on that set
to an order preserving homeomorphism of $\I$, first by using
continuity to extend to the closure of the embedded set and then using
affine extensions on the complementary intervals of this closure.

The approach we take is somewhat similar.  We 
consider the group $\Z^n$ of $n$-tuples of integers and provide it with
a linear order $\succ$ which is the lexicographic ordering, i.e. \
$(x_1, \dots, x_n) \prec (y_1, \dots, y_n)$ if and only if $x_i = y_i$
for $ 1 \le i < k$ and $x_k < y_k$ for some $0 \le k \le n.$ 
It is well known and easy to show that the standard linear 
action of $\N_n$ on $\Z^n$ preserves this ordering.

Instead of embedding $\Z^n$ as a countable set of points, however,
for each $(q_1, \dots, q_n) \in \Z^n$ we will embed a closed
interval $I(q_1, \dots, q_n)$ in $\I$.  We do this in such a way that the 
intervals are disjoint except that each one intersects its successor
in a common endpoint, and so that the order of the intervals
in $\I$ agrees with the lexicographic ordering on $\Z^n.$   We also
arrange that the complement of the union of these intervals is
a countable set.  If for each $\alpha \in \N_n$ we defined $g_{\alpha}$
on $I(q_1, \dots, q_n)$ to be the canonical affine map to
$I(\alpha(q_1, \dots, q_n))$  this would define $g_{\alpha}$ on
a dense subset of $I$ and it would extend uniquely to a homeomorphism,
giving an action of $\N_n$ by homeomorphisms.

In order to improve this to a $C^1$ action we must replace the affine
maps from $I(q_1, \dots, q_n)$ to $I(\alpha(q_1, \dots, q_n))$ with
elements of another canonical family of diffeomorphisms which has two
key properties. The first is that these subinterval diffeomorphisms
must have derivative $1$ at both endpoints in order to fit together in
a $C^1$ map.  The second is that if $I_k$ is a strictly monotonic
sequence of these subintervals then the restriction of $g_{\alpha}$ to
$I_k$ must have a derivative which converges uniformly to $1$ as $k$
tends to infinity.  In particular, by the mean value theorem the ratio
of the lengths of $\alpha(I_k)$ and $I_k$ must tend to $1$ for any
element $\alpha \in \N_n$.  This makes the choice of the lengths of
these intervals one of the key ingredients of the proof.  We choose
these lengths with a positive parameter $K$.  We are then able to show
that as this parameter increases to infinity the $C^1$ size of any
$g_{\alpha}$ goes to zero.  This allows us to prove the existence of
an effective $C^1$ action of $\N_n$ on $\I$ with generators chosen
from an arbitrary $C^1$ neighborhood of the identity.

We proceed to the details beginning with some calculus lemmas.  The
first of them develops a family of interval diffeomorphisms which we
will us as a replacement for affine interval maps as we discussed
above.  We are indebted to Jean-Christophe Yoccoz for the proof of the
following lemma which substantially simplifies our earlier approach.

\begin{lem}\label{lem:phi}
For each $a$ and $b \in (0, \infty)$ there 
exists a  $C^1$ orientation preserving 
diffeomorphism $\phi_{a,b}: [0,a] \to [0,b]$ with the
following properties:
\begin{enumerate}
\item For any $a,b,c \in (0, \infty)$ and for any 
$x \in [0,a],\ \phi_{b,c}(\phi_{a,b}(x)) = \phi_{a,c}(x).$
\item For all $a, b, \ \phi'_{a,b}(0) = \phi'_{a,b}(1) = 1.$
\item Given $\varepsilon >0$ there exists $\delta >0$ such that 
for all $x \in [0,a],$
\[
\Big |\phi'_{a,b}(x) - 1 \Big | < \varepsilon, \text{ whenever }
\Big |\frac{b}{a} - 1 \Big | < \delta
\]
\end{enumerate}
\end{lem}

\begin{proof}
For $u_0 \in (0, \infty)$  we define 
\begin{equation*}
L(u_0) = \int_{-\infty}^{\infty} \frac{du}{u^2 + u_0^2}.
\end{equation*}
A change of variables $u = u_0 x$ allows one to conclude that
\begin{equation*}
L(v) = \frac{1}{u_0}\int_{-\infty}^{\infty} 
\frac{dx}{x^2 + 1} = \frac{\pi}{u_0}
\end{equation*}

For each $u_0 \in (0,\infty)$ we define $\psi_{u_0}: \R \to (0,L(u_0))$ by
\begin{equation*}
\psi_{u_0}( t) = \int_{-\infty}^{t} \frac{du}{u^2 + u_0^2}
\end{equation*}
Clearly $\psi_{u_0}$ is a real analytic diffeomorphism from $\R$ onto
the open interval $(0,L(u_0)).$

Now define $\phi_{a,b}:[0,a] \to [0,b]$ by setting 
$\phi_{a,b}(0) = 0, \phi_{a,b}(a) = b$ and for $x \in (a,b)$ letting
$\phi_{a,b}(x) = \psi_{u_1}(\psi_{u_0}^{-1}(x))$ where
$a = L(u_0)$ and $b = L(u_1).$
This is clearly
a homeomorphism of $[0,a]$ to $[0,b]$ which is real analytic
on $(0,a).$  We observe that the derivative of $\phi_{a,b}(x)$ at
the point $x = \psi_{u_0}(t)$ is given by
\[
\phi_{a,b}'(x) = \frac{\psi_{u_1}'(t)}{\psi_{u_0}'(t)} =
\frac{t^2 + u_0^2}{t^2 + u_1^2}
\]

{From} this it is clear that $\phi_{a,b}'(x)$ can be continuously
extended to the endpoints of the interval $[0,a]$ by assigning
it the value $1$ there.  Hence $\phi_{a,b}(x)$ is $C^1$ on the
closed interval $[0,a]$ and satisfies property (2) above.
Property (1) is clear from the definition.

To check property (3) we observe that if, as above $x = \psi_{u_0}(t)$
then
\[
\Big |\phi'_{a,b}(x) - 1 \Big | = \Big |
\frac{t^2 + u_0^2}{t^2 + u_1^2} - 1 \Big |
\]
which is easily seen to assume its maximum when $t = 0.$  But
\[
\Big |\phi'_{a,b}(x) - 1 \Big | \le \Big |\frac{u_0^2}{u_1^2} - 1 \Big |
=\Big |\frac{b^2}{a^2} - 1 \Big |
\]
since $u_0 = \pi/a$ and $u_1 = \pi/b.$   From this it is clear that
property (3) holds.
\end{proof}

We will also need the following technical lemma.

\begin{lem}\label{lem:xy.lim}
Suppose $n$ is a positive even integer and $K > 0.$
Then for $(x,y) \in \R^2$
\[
\lim_{\|(x,y)\| \to \infty} \frac{|(x + y)^n - y^n|}{x^{n+2} + y^{n} + K} = 0.
\]
Moreover, given $\varepsilon >0$, for $K$ sufficiently large 
\[
\frac{|(x + y)^n - y^n|}{x^{n+2} + y^{n} + K} < \varepsilon
\]
for all $(x,y) \in \R^2.$
\end{lem}
\begin{proof}
The numerator $(x + y)^n - y^n$ is a sum of mononomials of the form
$C x^k y^{n-k}$ where $0 < k \le n$.  Hence it suffices to prove
\[
\lim_{\|(x,y)\| \to \infty} \frac{|x|^k|y|^{n-k}}{x^{n+2} + y^{n} + K} = 0
\]
for $0 < k \le n$.

If for some $\varepsilon >0$ 
\[
\frac{|x|^k|y|^{n-k}}{x^{n+2} + y^{n} + K} \ge \varepsilon
\]
we want to  show that there is an upper bound for $|x|$ and $|y|$.

We first observe that
\begin{equation*}
\frac{|x|^k}{|y|^{k}} =  \frac{|x|^k|y|^{n-k}}{|y|^n}
> \frac{|x|^k|y|^{n-k}}{x^{n+2} + y^{n} + K}  \ge \varepsilon  
\end{equation*}
or
\begin{equation}
|x| > \varepsilon^{1/k} |y| 
\label{x>y}
\end{equation}

Similarly,
\begin{equation*}
\frac{|y|^{n-k}}{|x|^{n-k +2}} = \frac{|x|^k|y|^{n-k}}{|x|^{n+2}} >
\frac{|x|^k|y|^{n-k}}{x^{n+2} + y^{n} + K} \ge \varepsilon
\end{equation*}
This implies that 
\begin{equation}
|y|^{n-k}> \varepsilon |x|^{n-k+2} \label{y>x.1}
\end{equation}
and we note that $n- k +2 > 0.$  In the case that $k = n$
we note that equation (\ref{y>x.1}) implies that $|x|$ is bounded and then
equation (\ref{x>y}) implies that $|y|$ is bounded.

On the other hand, when $k < n$ we have 
\begin{equation}
|y| > E |x|^d 
\label{y>x.2}
\end{equation}
where 
\begin{equation*}
E = \varepsilon^{\frac{1}{n-k}} \text{ and } 
d = \frac{n -k +2}{n -k} = 1 + \frac{2}{n-k}.
\end{equation*}

Combining equations (\ref{x>y}) and (\ref{y>x.2}) we see 
\begin{equation*}
|x| > \varepsilon^{1/k} |y| > \varepsilon^{1/k} E |x|^{d}
\end{equation*}

Since $d >1$ this clearly implies that $|x|$ is
bounded by a constant depending only on $n$ and $\varepsilon.$ Then
equation (\ref{x>y}) implies that $|y|$ is also bounded.  The
contrapositive of these assertions is that for $\|(x,y)\|$
sufficiently large we have
\begin{equation*}
\frac{|x|^k|y|^{n-k}}{x^{n+2} + y^{n} + K} < \varepsilon
\end{equation*}
which implies the desired limit.

To prove the second assertion of the lemma we observe that we have shown
there is a constant $M>0$ such that
\[
\frac{(x + y)^n - y^n}{x^{n+2} + y^{n} + K} < \varepsilon
\]
whenever $\|(x,y)\| >M.$  Morever the constant $M$ is independent of $K$.
Hence clearly if $K$ is sufficiently large, this inequality
will hold for all $(x,y).$
\end{proof}

Consider the group $\Z^n$ of $n$-tuples of integers and provide it with
a linear order $\succ$ which is the lexicographic ordering, i.e. \
$(x_1, \dots, x_n) \prec (y_1, \dots, y_n)$ if and only if $x_i = y_i$
for $ 1 \le i < k$ and $x_k < y_k$ for some $0 \le k \le n.$ 
We are now prepared to define a set of closed intervals in one-to-one
correspondence with elements of $\Z^n$.  As mentioned above we will
make this definition with a positive parameter $K$ which will be 
used subsequently to get generators of our action in an arbitrary
$C^1$ neighborhood of the identity.

\begin{defn}

For $K >0$  let $B_K : \Z^{n} \to \R$ be defined by
\begin{align*}
B_K( q_1, q_2, \dots, q_n) &= K +  \sum_{j=1}^{n} q_j^{4n-2j+2}\\
&= q_1^{4n+2} + q_2^{4n} + \dots + q_{n-1}^{2n+4} + q_n^{2n+2} + K.
\end{align*}

We will show convergence of the series

\[
S_K = \sum_{( q_1, q_2, \dots, q_n) \in \Z^n} \frac{1}{B_K( q_1, q_2,
\dots, q_n)}.
\]
For $( r_1, r_2, \dots, r_n)$ we define $S_K( r_1, r_2, \dots, r_n)$ by
\[
S_K( r_1, r_2, \dots, r_n) = \sum_{( q_1, q_2, \dots, q_n)
\prec ( r_1, r_2, \dots, r_n)} \frac{1}{B_K( q_1, q_2, \dots, q_n)}.
\]
Finally we define the closed interval 
\[
I_K( r_1, r_2, \dots, r_n) = [S_K( r_1, r_2, \dots, r_n),
S_K( r_1, r_2, \dots, r_n + 1)].
\]
\end{defn}

\begin{remark}
In an earlier version of this paper we used $2n-2j +2$ as the
the exponent in the definition of $B_K$ instead of $4n-2j +2$.
As a result (as was pointed out to us) the series $S_K$ might 
not converge.  However, the only properties we use of $B_K$ are that
$S_K$ converges and that the ratio of the value of $B_K$ at certain points
limits to $1$.  These limits are shown in Lemma~\ref{lem:B.lim} and Lemma~\ref{lem:B.est}.
\end{remark}

We are indebted to Nir Avni and Elton Hsu for the following argument 
showing the  series $S_K$ converges.   We will use comparison with 
a convergent integral.
So we need to show
\[
\int_{\R^n} \frac{1}{B_K(q)} dV < \infty.
\]
where $dV$ is the standard volume element on $\R^n$ and $r(q) = \|q\|$ is the
standard norm.  
Clearly it suffices to show finiteness of the integral over the region $\|q\| \ge 1$.
Define $w_i(q) = q_i / r(q)$ and let
\[
C_n(q) = \sum_{i=1}^n w_i(q)^{4n-2i +2}.
\]
Note that $\sum w_i^2 = 1$, so for each $q$ there is some $j$ 
with $|w_j(q)| \ge 1/\sqrt{n}$ and hence  $C_n$ is bounded below
(with a bound depending on $n$).

We now observe that if $r \ge 1$
\begin{align*}
B_K( q_1, q_2, \dots, q_n) - K  &= \sum_{i=1}^{n} q_i^{4n-2i+2}\\
&= \sum_{i=1}^{n} r^{4n-2i+2} w_i^{4n-2i+2}\\
&\ge r^{2n+2}   \sum_{i=1}^{n}w_i^{4n-2i+2}\\
&\ge C_n(q) r(q)^{2n+2}.
\end{align*}

Let $C$ be a lower bound for $C_n(q)$.
Then since $dV = r^{n-1}\ dr \wedge \omega$ where $\omega$
is a volume form on the sphere $\|q\| = 1$, we have 
\begin{align*}
\int_{\|q\| \ge 1} \frac{1}{B_K(q)} dV &\le \int_{\|q\| \ge 1} \frac{1}{K + C_n(q) r^{2n+2}(q)}dV\\
&\le \int_{\|q\| \ge 1} \frac{r^{n-1}\ dr \wedge \omega}{K +C r^{2n+2}} \\
&\le \int_{\|q\| \ge 1} \frac{1}{C r^{n+3}} dr \wedge \omega < \infty\\
&= \frac{1}{C} \int_1^\infty  \frac{1}{ r^{n+3}}\ dr  \int_{\|q\| = 1} \omega   < \infty.
\end{align*}
It follows by  comparison that the series $S_K$ converges.

We will
denote the interval $[0,S_K]$ by $I_K$ and we note that it is nearly
the union of the intervals $I_K( r_1, r_2, \dots, r_n)$.

\begin{lem}
\label{lem:def.J}
There is a countable closed set $J_K$ in the interval $I_K=[0,S_K]$,
such that 
\[
I_K = J_K \cup \Big( \bigcup_{( q_1, q_2, \dots, q_n) \in \Z^n}
I_K( q_1, q_2, \dots, q_n)\Big)
\]
\end{lem}
\begin{proof}
Define $J_K$ by the equality in the statement of the lemma.  From the
definitions it follows that every point $J_k$ is the limit of one of the 
Cauchy sequences $\{S_K(r_1,\ldots ,r_i,\ldots r_n)\}$ as $r_i$ ranges
from $1$ to either $\infty$ or $-\infty$ and the other $r_j, j\neq i$
are fixed.  It follows that $J_K$ is countable.  The set $J_K$ is closed 
since its complement is clearly open.
\end{proof}

We note that the intervals $\{I_K( q_1, q_2, \dots, q_n)\}$ occur in 
the interval $I_K$ in exactly the order given by $\prec.$  They
have disjoint interiors but each such interval intersects its successor
in a common endpoint.

\begin{defn}
Let $\nu : I_K \setminus J_K \to \Z^{n}$ be defined by setting 
$\nu(x) = ( q_1, q_2, \dots, q_n)$ if $x \in I_K(  q_1, q_2, \dots, q_n)$ and
$x \notin I_K(  q_1, q_2, \dots, q_{n +1}).$
\end{defn}

\begin{lem} \label{lem:def.g_i}
If $K$ is sufficiently large, then for every 
$\alpha \in \N_n$ there exists a homeomorphism
$g_\alpha : I_K \to I_K$ such that for each
$(q_1, \dots, q_n)\in \Z^n,$ 
\begin{enumerate}
\item $g_\alpha(I_K(q_1, \dots, q_n)) = I_K(\alpha(q_1, \dots, q_n)),$ and
\item $g_{\alpha \beta}(x) = g_{\alpha}(g_{\beta}(x))$ for all $\alpha, \beta
\in \N_n$ and all $x \in I_K.$
\item for $x \in I_K(q_1, \dots, q_n),$
$g_\alpha(x) = \phi_{b_1,b_2}( x - c_1) + c_2$, 
where 
\begin{align*}
b_1 &= |I_K(q_1, \dots, q_n)| = B_K(q_1, \dots, q_n)^{-1},\\
b_2 &= |I_K(\alpha(q_1, \dots, q_n)| = B_K(\alpha(q_1, \dots, q_n))^{-1},\\
c_1 &= S_K(q_1, \dots, q_n), \text{ and }\\
c_2 &= S_K(\alpha(q_1, \dots, q_n)).
\end{align*}
and $\phi_{b_1,b_2}: [0,b_1] \to [0,b_2]$ is the diffeomorphism guaranteed
by Lemma \ref{lem:phi}.
\end{enumerate}
\end{lem}
\begin{proof}
For $x \in I_K(q_1, \dots, q_n),$ we define
$g_\alpha(x)$ to be $\phi_{b_2}( \phi_{b_1}^{-1}(x - c_1)) + c_2$.
The values of the constants $b_j, c_j,$ for $j = 1,2$ were
chosen so that (1) holds. 
This gives a homeomorphism $g_\alpha$ from $I_K \setminus J_K$ to itself
which preserves the order on the interval $I_K$.  Since $J_K$ is a
closed countable subset of $I_K$ we may extend $g_\alpha$ to all of 
$I_K$ in the unique way which makes it an order preserving
function on $I_K$.  Clearly this makes $g_\alpha$ a homeomorphism
of $I_K$ to itself.

In order to show property (2) we let 
$b_3 = B_K({\beta \alpha}(q_1, \dots, q_n))^{-1}$ and
$c_3 = S_K({\beta \alpha}(q_1, \dots, q_n)).$  Then
if $x \in I_K(q_1, \dots, q_n)$,
\begin{align*}
g_{\beta}(g_{\alpha}(x)) &= \phi_{b_2,b_3}(g_{\alpha}(x) - c_2) + c_3\\
&= \phi_{b_2,b_3}( \phi_{b_1,b_2}( (x - c_1) + c_2  - c_2) + c_3,
\text{ by Lemma \ref{lem:phi}}\\
&= \phi_{b_1,b_3}( x - c_1)  + c_3\\
&= g_{\beta \alpha}(x).
\end{align*}

Since this holds for any $(q_1, \dots, q_n)\in \Z^n$ we have shown
$g_{\beta}(g_{\alpha}(x)) = g_{\beta \alpha}(x)$ for a dense set of
$x \in I_K$.  Continuity then implies this holds for all
$x \in I_K$.
\end{proof}

\begin{defn}
For $1 \le i <n$ let $\sigma_i \in \N_n$ denote the matrix
with all entries on the diagonal equal to $1$, with entry
$(i+1, i)$ equal to $1$, and with all other entries $0.$
We will denote the homeomorphism $g_{\sigma_i}: I_K \to I_K$ by $g_i$.
\end{defn}

The elements $\{\sigma_i\}$  form a set of generators for the group $\N_n$.  Our
next objective is to show that in fact $g_i = g_{\sigma_i}$ is a $C^1$
diffeomorphism.  For this we will need a sequence of technical lemmas.

\begin{lem}
\label{lem:B.lim}
Suppose ${x_k}$ is a monotonic sequence in $I_K \setminus J_K$ converging to a
point of $J_K.$  Then for each $1 \le i < n$
\[
\lim_{k \to \infty} \frac{B_K(\sigma_i(\nu(x_k)))}{B_K(\nu(x_k))} = 1.
\]
\end{lem}
\begin{proof}
We may assume without loss of generality that the sequence $\{x_k\}$
is monotonic increasing.  For $1 \le j \le n, k > 0$ we define
$q_j(k)$ by $\nu(x_k) = ( q_1(k) , q_2(k) , \dots, q_n(k))$.  Then
each sequence $\{q_j(k)\}_{k=1}^\infty$ is monotonic increasing.  At
least one of these sequences is unbounded, since any bounded ones are
eventually constant and if $\{q_j(k)\}$ were eventually constant for
all $j \le n$ then the sequence of intervals $I_K( \nu( x_k))$ would
be eventually constant and the limit of the sequence $\{x_k\}$ would
be in the final interval and hence not in $J_K$. Let $r$ be the
smallest $j$ such that $\{q_j(k)\}$ is unbounded.  Then since
$\{x_k\}$ is monotonic increasing we have $\lim_{k \to \infty} q_r(k)
= \infty$ and for $j < r$ the sequence $\{q_j(k)\}$ is eventually
constant.

We note that
\[
B_K(\sigma_i(\nu(x_k))) = 
B_K( q_1(k),\dots, q_{i}(k), q_{i+1}(k) + q_{i}(k), q_{i+2}(k), \dots,  q_n(k))
\]
and hence that 
\begin{align*}
\frac{B_K(\nu(x_k)) - B_K(\sigma_i(\nu(x_k)))}{B_K(\nu(x_k))}
&= \frac{(q_{i+1}(k) + q_{i}(k))^{4n -2i} - q_{i+1}(k)^{4n -2i}}
{B_K(q_1(k), \dots  q_n(k))} \\
&=\frac{ P(q_{i}(k),  q_{i+1}(k))}{B_K(q_1(k), \dots  q_n(k))}
\end{align*}
where $P(x,y) = (x + y)^{4n-2i} - y^{4n-2i}$.

We first consider the case that $r > i+1.$  Then for large $k$ we have that
$P(q_i(k),  q_{i+1}(k))$ is bounded (in fact eventually constant) so 
\[
\lim_{k \to \infty} \frac{B_K(\nu(x_k)) - B_K(\sigma_i(\nu(x_k)))}{B_K(\nu(x_k))}
= \lim_{k \to \infty} \frac{P(q_i(k), q_{i+1}(k))} {B_K(
q_1(k),\dots, q_n(k))} = 0
\]
and we have the desired result.

In case $r = i+1$ or $r = i$ we observe
\[
\frac{P(q_i(k),  q_{i+1}(k))}
{ q_1(k)^{4n} + \dots + q_n(k)^{2n} + K} \le
\frac{P(q_i(k),  q_{i+1}(k))}
{ q_i(k)^{4n-2i+2} + q_{i+1}(k)^{4n-2i} + K}
\]

and at least one of $q_i(k)^{4n-2i +2}$ and $q_{i+1}(k)^{4n-2i}$
tends to infinity so by Lemma \ref{lem:xy.lim} we again have
\[
\lim_{k \to \infty} \frac{B_K(\nu(x_k)) -
B_K(\sigma_i(\nu(x_k)))}
{B_K(\nu(x_k))}
=\lim_{k \to \infty} \frac{P(q_i(k),  q_{i+1}(k))}
{B_K( q_1(k),\dots,  q_n(k))} = 0
\]

Finally in case $r<i$ 
\[
\frac{ |P(q_i(k),  q_{i+1}(k))|}
{B_K( q_1(k),\dots,  q_n(k))} \le
\frac{ |P(q_i(k),  q_{i+1}(k))|}
{q_r(k)^{4n-2r+2} +q_i(k)^{4n-2i+2} + q_{i+1}(k)^{4n-2i} + K}
\]
and $|q_r(k)|$ tends to infinity.  But given $\varepsilon >0$ by
Lemma \ref{lem:xy.lim}  there is an $M>0$ such that whenever 
$\|(q_i(k),  q_{i+1}(k))\| >M$ we have 
\begin{eqnarray*}
\frac{ |P(q_i(k),  q_{i+1}(k))|}
{q_r(k)^{4n-2r+2} +q_i(k)^{4n-2i +2} + q_{i+1}(k)^{4n-2i} + K}\\
\le \frac{ |P(q_i(k),  q_{i+1}(k))|}{q_i(k)^{4n-2i +2} + 
q_{i+1}(k)^{4n-2i} + K} < \varepsilon
\end{eqnarray*}
On the other hand if $\|(q_i(k), q_{i+1}(k))\| \le M$ and
$|q_r(k)|$ is sufficiently large
\[
\frac{ |P(q_i(k),  q_{i+1}(k))|}
{q_r(k)^{4n-2r+2} +q_i(k)^{4n-2i+2} + q_{i+1}(k)^{4n-2i} + K} < \varepsilon
\]
So in all cases we have the desired limit
\begin{equation}
\lim_{k \to \infty} \frac{|B_K(\nu(x_k)) -
B_K(\sigma_i(\nu(x_k)))|}
{B_K(\nu(x_k))} = 0 
\label{B.lim.1}
\end{equation}
which implies
\begin{equation*}
\lim_{k \to \infty} \Big |1 - \frac{B_K(\sigma_i(\nu(x_k)))}{B_K(\nu(x_k))}\Big |
= 0
\end{equation*}
and hence that
\[
\lim_{k \to \infty} \frac{B_K(\sigma_i(\nu(x_k)))}{B_K(\nu(x_k))} = 1.
\]
\end{proof}

\begin{lem}
\label{lem:deriv.lim}
Suppose ${x_k}$ is a monotonic sequence in $I_K \setminus J_K$ converging to a
point of $J_K.$  Then for each $1 \le i < n,$ 
\[
\lim_{k \to \infty} g_i'(x_k) = 1.
\]
\end{lem}
\begin{proof}

For $x \in I_K(q_1, \dots, q_n),\ g_i$ is defined by
$g_i(x) = \phi_{b_1,b_2}(x - c_2) + c_1$, for
some constants $c_1$ and $c_2$,
where 
\[
b_1 = b_1(k) = \frac{1}{B_K(\nu(x_k))} \text{ and } b_2 = b_2(k) = 
\frac{1}{B_K(\sigma_i(\nu(x_k)))}
\]
So $g_i'(x_k) = \phi_{b_1,b_2}'(x_k)$.

We wish to show that given $\varepsilon >0$ for $k$ sufficiently large
$|1 - g_i'(x_k)| < \varepsilon$.
To do this, by Lemma \ref{lem:phi} we need only show that for any $\delta >0$ 
\[
\Big |\frac{b_2(k)}{b_1(k)} - 1 \Big | < \delta 
\]
for sufficiently large $k.$

But
\begin{equation*}
\Big |\frac{b_2(k)}{b_1(k)} - 1 \Big |
= \Big |\frac{B_K(\nu(x_k))}{B_K(\sigma_i(\nu(x_k)))} - 1 \Big |
\end{equation*}
which, by Lemma \ref{lem:B.lim} tends to $0$ as $k$ tends to infinity.
Hence
\[
\Big |\frac{b_2(k)}{b_1(k)} - 1 \Big | < \delta
\]
for $k$ sufficiently large.
\end{proof}

The proofs of Lemmas (\ref{lem:B.est}) and (\ref{lem:deriv.est})
below are parallel to those of Lemmas (\ref{lem:B.lim}) and
(\ref{lem:deriv.lim}).  The difference is that before we were interested in
the limiting values of some quantities as a sequence 
$\{x_k\}$ of values of $x$ in $I_K \setminus J_K$ converged, while now we
are interested in estimating the same quantities but, uniformly in
$x$, as the parameter $K$ tends to infinity.

\begin{lem}
\label{lem:B.est}
Suppose $1 \le i < n$.
Given $\varepsilon >0$ there is a $K_0 >0$ such that whenever $K > K_0$
\[
\Big |\frac{B_K(\sigma_i(\nu(x)))}{B_K(\nu(x))} - 1 \Big | < \varepsilon
\]
for all $x \in I_K \setminus J_K.$
\end{lem}
\begin{proof}
For $1 \le j \le n,$ we define $q_j(x)$ by $\nu(x) = ( q_1(x) ,
q_2(x) , \dots, q_n(x))$.

We note that
\[
B_K(\sigma_i(\nu(x))) = 
B_K( q_1(x),\dots, q_{i}(x), q_{i+1}(x)+ q_{i}(x), q_{i+2}(x), \dots,  q_n(x))
\]
and hence if  $P(u,v) = (u + v)^{4n-2i} - v^{4n-2i}$, we have
\begin{align*}
\Big |\frac{B_K(\nu(x)) - B_K(\sigma_i(\nu(x)))}{B_K(\nu(x)}\Big |
&= \Big |\frac{(q_{i+1}(x) + q_{i}(x))^{4n -2i} - q_{i+1}(x)^{4n -2i}}
{B_K(q_1(x), \dots  q_n(x)}\Big | \\
&=\frac{ |P(q_i(x),  q_{i+1}(x)|}
{B_K(q_1(x), \dots  q_n(x))}\\
&\le\frac{ |P(q_i(x),  q_{i+1}(x))|}
{q_{i}(x)^{4n -2i +2} + q_{i+1}(x)^{4n -2i}+ K}.
\end{align*}
So by the second part of Lemma (\ref{lem:xy.lim}) we conclude that
if $K$ is sufficiently large
\begin{equation}
\Big |\frac{B_K(\nu(x)) - B_K(\sigma_i(\nu(x)))}{B_K(\nu(x))}\Big |
< \varepsilon. \label{B.est.1}
\end{equation}

{From} this we see
\begin{equation*}
\Big | 1 - \frac{B_K(\sigma_i(\nu(x)))}{B_K(\nu(x))}\Big | < 
\varepsilon
\end{equation*}
as desired.
\end{proof}

\begin{lem}
\label{lem:deriv.est}
Given $\varepsilon >0,$ if $K$ is chosen sufficiently large then for
all $x \in I_K \setminus J_K$ and all $1 \le i < n$
\[
|g_i'(x) - 1| < \varepsilon.
\]
\end{lem}
\begin{proof}

For $x \in I_K(q_1, \dots, q_n),\ g_i$ is defined by
$g_i(x) = \phi_{b_1,b_2}(x - c_2) + c_1$, for
some constants $c_1$ and $c_2$,
where 
\[
b_1 = b_1(x) = \frac{1}{B_K(\nu(x))} \text{ and } b_2 = b_2(x) =
\frac{1}{B_K(\sigma_i(\nu(x)))}.
\]
So $g_i'(x) = \phi_{b_1,b_2}'(x).$

We wish to show that given $\varepsilon >0$, if $K$ is sufficiently large then
$|1 - g_i'(x)| < \varepsilon$.
To do this, by Lemma \ref{lem:phi} we need only show that given any $\delta >0$ 
there is a $K_0 >0$ such that $K > K_0$ implies
\[
\Big |\frac{b_2(x)}{b_1(x)} - 1 \Big | < \delta
\]
for all $x \in I_K \setminus J_K.$

But
\begin{equation*}
\Big |\frac{b_2(x)}{b_1(x)} - 1 \Big |
= \Big |\frac{B_K(\sigma_i(\nu(x)))}{B_K(\nu(x))} - 1\Big |
\end{equation*}
and the result follows from Lemma \ref{lem:B.est}.
\end{proof}

\begin{prop} \label{prop:gi.C1}
For $K$ sufficiently large the homeomorphism $g_i : I_K \to I_K$ is a
$C^1$ diffeomorphism with derivative $1$ at both endpoints.  Given
$\varepsilon >0$ there exists $K_0$ such that whenever $K > K_0$ we
have
\[
|g_i'(x) - 1| < \varepsilon
\]
for all $x \in I_K.$
\end{prop}
\begin{proof}
We know that the function $f(x) = g_i'(x)$ exists and is continuous
on $I_K \setminus J_K.$  By Lemma (\ref{lem:deriv.lim}) we can extend
it continuously to all of $I_K$ by setting $f(x) = 1$ for $x \in J_K.$

We define a $C^1$ function $F$ by
\[
F(x) = \int_0^x f(t) dt
\]
and will show that $F(x) = g_i(x).$  To see this let
$\phi(x) = g_i(x) - F(x).$  Then $\phi(0) = 0$ and $\phi(x)$ is a continuous
function whose derivative exists and is $0$ on $I_K \setminus J_K.$
Since $J_K$ is countable $\phi( J_K)$ has measure zero.  But 
$I_K \setminus J_K$ has countably many components on each of which
$\phi$ is constant.  It follows that $\phi(I_K)$ has measure zero
and hence $\phi(I_K) = \{0\}.$  Therefore $g_i(x) = F(x)$  is $C^1.$

Lemma (\ref{lem:deriv.est}) and the fact that
$g_i'(x) = 1$ for $x \in J_K$ imply that for $K$ sufficiently large
\[
|g_i'(x) - 1| < \varepsilon
\]
for all $x \in I_K.$   The inverse function theorem then implies that
$g_i$ is a $C^1$ diffeomorphism.
\end{proof}

Recall that we have given the group $\Z^n$ the 
lexicographic ordering,  $\succ$ i.e. \
$(x_1, \dots, x_n) \succ (y_1, \dots, y_n)$ if and only if $x_i = y_i$
for $ 1 \le i < k$ and $x_k > y_k$ for some $0 \le k \le n.$ We note
that this order is translation invariant, indeed $(x_1, \dots, x_n)
\succ (y_1, \dots, y_n)$ if and only if $(x_1-y_1, \dots, x_n-y_n) \succ
(0, \dots, 0).$ 

It is well known and easy to check that the nilpotent
group $\N_n$ of $n\times n$ lower-triangular integer matrices with ones
on the diagonal acts on $\Z^n$ preserving the order $\succ$.

\begin{thm}
\label{theorem:nilpotent.C1}
The group $\N_n$ is isomorphic to a subgroup of $\Diff^1(M)$ for
$M=\R,\T^1,$ or $\I$.  For $M=\T^1$ or $\I$ the elements of this
subgroup corresponding to the generators $\{\sigma_i\}$ of $\N_n$ may
be chosen to be in an arbitrary neighborhood of the identity in
$\Diff^1(M)$.
\end{thm}

\begin{proof}
We first consider the case that $M = \I.$  Given $\varepsilon > 0$ we
choose $K$ sufficiently large that the conclusion of
Proposition (\ref{prop:gi.C1}) holds.  We define 
$\Phi: \N_n \to \Homeo(I_K)$ by $\Phi(\alpha) = g_\alpha.$
Lemma (\ref{lem:def.g_i}) asserts that $\Phi$ is an injective homomorphism.
Proposition (\ref{prop:gi.C1}) asserts that $g_i = \Phi(\sigma_1)$ is a
$C^1$ diffeomorphism so in fact $\Phi(\N_n)$ lies in $\Diff^1(I_K)$.

Define the injective homomorphism $\Psi: \N_n \to \Diff^1(\I)$ by 
$\Psi(\alpha)(x) = S_K \Phi(\alpha)(x/S_K) = S_K g_\alpha(x/S_K).$
Restricting this action to the interior of $\I$ givens an action
on $\R.$ 

Note that every element of $\Psi(\N_n)$ has derivative $1$ at both
endpoints of $\I$.  Gluing endpoints together gives an action on $S^1$.

In the case $M=\T^1$ or $\I$ we clearly have $|\Psi(\sigma_i)'(x) - 1|
< \varepsilon$ by Proposition (\ref{prop:gi.C1})
\end{proof}

Now every finitely-generated, torsion-free nilpotent group $N$ is
isomorphic to a subgroup of $\N_n$ for some $n$ (see Theorem 4.12 of
\cite{Ra} and its proof).  This together with Theorem 
\ref{theorem:nilpotent.C1}
immediately implies Theorem \ref{thm:fg.nilpotent}.

\subsection{$C^\infty$ actions on $\R$}

For certain nilpotent groups we can give a $C^\infty$ action on $\R.$
But we will see in the next section that for $n > 2$ there is no
$C^2$ action of $\N_n$ on $\R.$

\bigskip
\noindent
{\bf Proof of Theorem \ref{theorem:smooth:example}: }
Choose a non-trivial $C^\infty$ diffeomorphism $\alpha$ of $[0,1]$ to
itself such that for $j = 0,1$ we have $\alpha(j) = j,\ \alpha'(j) =
1,$ and $\alpha^{[k]}(j) = 0,$ for all $k >1.$

Define three $C^\infty$ diffeomorphisms $f,h_0,h_1: \R \to \R$ by
$$
\begin{array}{ll}
f(x) = x - 1& \\
h_0(x) = \alpha(x - m) + m &\mbox{\ for\ }x \in [m, m+1]\\
h_1(x) = \alpha^{m}(x - m) + m&\mbox{\ for\ }x \in [m, m+1]
\end{array}
$$ 
Then it is easy
to check that $h_0$ and $h_1$ commute, $f$ and $h_0$ commute, and
$[f,h_1] = f^{-1}h_1^{-1}fh_1 = h_0$.  Hence the group generated by
$f$ and $h_1$ is nilpotent with degree of nilpotency $2$.

Given any $n >0$, we will inductively define $h_k$ for $0 \le k \le n$ in
such a way that $h(m) = m$ for $m \in \Z$, $[f,h_k] = h_{k-1}$ for 
$k>1,$ and $h_i h_j = h_j h_i.$ We do this by letting
$h_k(x) = x$ for $x \in [0, 1]$ and recursively defining $h_k^{-1}$.
We first define it for $x > 1$ by
\[
h_k^{-1}(x) =  h_{k-1}^{-1}f^{-1}h_k^{-1}f(x) \text { for } x > 1.
\]
Note this is well defined recursively  because 
if $x \in [n,n+1]$ the right hand side requires only that
we know the value of $h_k^{-1}(f(x))$ and $f(x) \in [n-1,n].$
We also note that $h_k^{-1} =  h_{k-1}^{-1}f^{-1}h_k^{-1}f$ implies
$h_{k-1}h_k^{-1} =  f^{-1}h_k^{-1}f$, so $h_{k-1}(x) =  
f^{-1}h_k^{-1}fh_k(x)$ for all $x > 1.$

Negative values of $x$ are handled similarly.  We define 
\[
h_k^{-1}(x) =  f h_{k-1}h_k^{-1}f^{-1}(x) \text { for }x < 0. 
\]
This is well defined recursively for $x<0$ because 
if $x \in [-n,-n+1]$ the right hand side requires only that
we know the value of $h_k^{-1}(f^{-1}(x))$ and $f^{-1}(x) \in 
[-n+1,-n+2].$ Again $h_k^{-1} =  f h_{k-1}h_k^{-1}f^{-1}$ implies
$f^{-1}h_k^{-1}f  =   h_{k-1}h_k^{-1}$ so $h_{k-1}(x) =  
f^{-1}h_k^{-1}fh_k(x)$ for all $x <0.$

We further note that from this definition one sees easily inductively 
that $h_k(x) = x$ for $x \in [-1,0]$ for all $k \ge 2.$
Finally we observe that this implies for $x \in [0,1]$ we
have $f^{-1}h_k^{-1}fh_k(x) =x$ so $f^{-1}h_k^{-1}fh_k(x) = h_{k-1}(x).$

It is clear from the relations $[f,h_k] = h_{k-1} \text{ for } k>1,\
[f,h_0] = id,$ and $h_i h_j = h_j h_i$ that the group generated by 
$f$ and the $h_k$ is nilpotent of degree at most $n$.  The degree of
nilpotency is at least $n$ since $h_0$ is a nontrivial 
$n$-fold commutator.
\endproof

If $n=2$, the group $G$ constructed above is isomorphic to the group
$\N_3$ of $3 \times 3$ lower triangular matrices.  In general, $G$ is
isomorphic to a semi-direct product of $Z^n$ and $Z$ where the $Z$
action on $Z^n$ is given by the $n \times n$ lower triangular matrix
with ones on the diagonal and subdiagonal and zeroes elsewhere.  It is
metabelian, i.e. solvable with derived length two.  We will see below
that this is a necessary condition for a smooth action on $\R.$

\subsection{Residually torsion-free nilpotent groups}
\label{section:residual}

In this section we consider some actions on $\R$ which are not irreducible,
but which are actions of a fairly large class groups.

\begin{defn}
A group $G$ is called {\em residually torsion-free nilpotent} if, for
every non-trivial element $g \in G$, there is a torsion-free nilpotent
group $N$ and a homomorphism $\phi: G \to N$ such that $\phi(g)$ is
non-trivial.  Equivalently, the intersection of every term in the lower
central series for $G$ is trivial.
\end{defn}

\bigskip
\noindent
{\bf Proof of Theorem \ref{theorem:residual}: }
Since $G$ is finitely generated it contains countably many elements.
Let $\{g_m \ | \ i\in \Z^+\}$ be an enumeration of the non-trivial 
elements of $G$. Let $\phi_m : G \to N_m$ be the homomorphism to
a torsion free nilpotent group guaranteed by residual nilpotence,
so $\phi_m(g_m)$ is non-trivial.  Replacing $N_m$ by $\phi_m(G)$
if necessary we may assume $N_m$ is finitely generated.

Let $I_m = [1/(m+1), 1/m]$.  Using Theorem \ref{thm:fg.nilpotent}, 
choose an effective action of $N_m$ by $C^1$ diffeomorphisms on the
interval $I_m.$ Recall that this action has the property that the
derivative of every element is $1$ at the endpoints of $I_m.$ In
addition, by Theorem \ref{thm:fg.nilpotent}, we may choose this action
so that the derivative of $\phi_m( g_i)$ satisfies
\[
| \phi_m( g_i)'(x) - 1 | < \frac{1}{2^m}
\]
for all $x \in I_m$ and all $1 \le i \le m.$

We then define the action of $G$ by $g(x) = \phi_m(g)(x)$ for
$x \in I_m$ and $g(0) = 0$ for all $g \in G.$
Clearly each $g_m$ is a $C^1$ diffeomorphism with $g_m'(0) = 1.$
The action is effective because $\phi_m(g)$ acts non-trivially on $I_m.$

Since every element of $G$ has derivative $1$ at both
endpoints of $\I$, we may glue the endpoints together 
to give a $C^1$ action on $\T^1$.
\endproof

It is an old result of Magnus that free groups and surface groups are
residually torsion-free nilpotent.   In particular Theorem
\ref{theorem:residual} gives a $C^1$ action of surface groups on $\R, \I$
and $\T^1$; we do not know another proof that such actions exist.  

Another interesting example is the {\em Torelli group} $T_g$, which is
defined to be the kernel of the natural action of the mapping class
group of a genus $g$ surface $\Sigma_g$ on $H_1(\Sigma_g,\Z)$.  For
$g\geq 3$, D. Johnson proved that $T_g$ is finitely-generated, and
Bass-Lubotzky proved that $T_g$ is residually torsion-free nilpotent.
Similarly, the kernel of the action of the outer automorphism group of a
free group is finitely-generated and residually torsion-free nilpotent.
Hence by Theorem \ref{theorem:residual} both of these groups are
subgroups of $\Diff^1(M)$ for $M=\R,\I,S^1$. In particular the Torelli
group is left orderable.  It is not known wether or not mapping class
groups are left-orederable, althought Thurston has proven that braid
groups are left-orderable.

\section{Restrictions on $C^2$ actions}

In the previous sections we showed that actions by nilpotent
subgroups of homeomorphisms are abundant. 
By way of contrast in the next two sections we show that
nilpotent groups of $C^2$ diffeomorphisms or PL homeomorphisms
are very restricted.  In this section we focus on $C^2$ actions.

\subsection{Kopell's Lemma}
Our primary tool is the following remarkable result of Nancy Kopell,
which is Lemma 1 of \cite{K}.

\begin{thm}[Kopell's Lemma] 
\label{thm:kopell}
Suppose $f$ and $g$ are $C^2$ orientation preserving
diffeomorphisms of an interval $[a,b)$ or $(a,b]$ and $fg = gf.$ If $f$ has no
fixed point in $(a,b)$ and $g$ has a fixed point in $(a,b)$ then
$g = id.$
\end{thm}

We will primarily use a consequence of this result which 
we now present.  

\begin{defn} 
We will denote by $\partial\Fix(f)$ the frontier of $\Fix(f)$,
i.e. the set $\partial\Fix(f) =\Fix(f) \setminus \Int(\Fix(f)).$
\end{defn}

\begin{lem}\label{lem:fix.component.R}
Suppose $f$ and $g$ are commuting orientation preserving $C^2$
diffeomorphisms of $\R$, each of which has a fixed point.  Then $f$
preserves each component of $\Fix(g)$ and vice versa. 
Moreover, $\partial\Fix(g) \subset \Fix(f)$
and vice versa.  The same result
is true for $C^2$ diffeomorphisms of a closed or half-open interval,
in which case the requirement that $f$ and $g$ have fixed points
is automatically satisfied by an endpoint.
\end{lem}

\begin{proof}
The proof is by contradiction.  Assume $X$ is a component of $\Fix(g)$
and $f(X) \ne X.$ Since $f$ and $g$ commute $f(\Fix(g)) = \Fix(g)$ so
$f(X)$ is some component of $\Fix(g)$.  Hence $f(X) \ne X$ implies $f(X)
\cap X = \emptyset.$ Let $x$ be an element of $X$ and without loss of
generality assume $f(x) < x.$ Define
\[
a = \lim_{n \to \infty} f^n(x) \text{ and } b = \lim_{n \to -\infty} f^n(x).
\]
Then at least one of the points $a$ and $b$ is finite since $f$ has a
fixed point.  Also $a$ and $b$ (if finite) are fixed under both $f$
and $g$, and $f$ has no fixed points in $(a,b).$ Thus Kopell's Lemma
\ref{thm:kopell} implies $g(y) = y$ for all $y \in (a,b)$
contradicting the hypothesis that $X$ is a component of $\Fix(g)$.
The observation that $\partial \Fix(f) \subset \Fix(g)$ follows from
the fact that components of $\Fix(f)$ are either points or closed
intervals so $x \in \partial \Fix(f)$ implies that either $\{x\}$ is
a component of $\Fix(f)$ or $x$ is the endpoint of an interval which
is a component of $\Fix(f)$ so, in either case, the fact that $g$
preserves this component implies $x \in \Fix(g).$

The proof in the cases that $f$ and $g$ are diffeomorphisms of a closed
or half-open interval is the same.
\end{proof}

Another useful tool is the following folklore theorem (see, e.g., 
\cite{FS} for a proof).

\begin{thm}[H\"{o}lder's Theorem] 
\label{thm:holder} 
Let $G$ be a group acting freely and 
effectively by homeomorphisms on any closed subset
of $\R$.  Then $G$ is abelian.
\end{thm}

\subsection{Measure and the translation number}

A basic property of finitely generated nilpotent groups of
homeomorphisms of a one-manifold $M$ is that they have invariant Borel
measures.  Of course in the case $M = \I$ or $T^1$ this follows from
the fact that nilpotent groups are amenable and any amenable group
acting on a compact Hausdorf space has an invariant Borel probability
measure.  In the case $M = \R$ this is a special case of a result due
to J. Plante \cite{P} who showed that any finitely generated subgroup
of $\Homeo(\R)$ with polynomial growth has an invariant measure which
is finite on compact sets.

We summarize these facts in the following 
\begin{theorem}
\label{thm:inv.measure}
Let $M=\R,\T^1,$ or $\I$.  Then any finitely generated nilpotent subgroup of
$\Homeo(M)$ has an invariant Borel measure $\mu$ which is finite
on compact sets.
\end{theorem}

If one has an invariant measure for a subgroup of $\Homeo(\R)$, it is
useful to consider the translation number which was discussed by 
J.~Plante in \cite{P}.  It is an analog of the rotation number for
circle homeomorphisms.

Suppose $G$ is a subgroup of $\Homeo(\R)$ which preserves a Borel
measure $\mu$ that is finite on compact sets.  Fix a point $x \in \R$
and for each $f \in G$ define
\begin{equation} \notag
\tau(f) =
\begin{cases}
\mu([x,f(x))	&\text{if $x < f(x)$} \\ 
0	&\text{if $x =f(x)$}\\
-\mu([f(x),x))	&\text{if $x >f(x)$} 
\end{cases}
\end{equation}

The function $\tau : G \to \R$ is called the {\em translation number}.
The following properties observed by J. Plante in \cite{P} are easy
to verify.
\begin{prop}\label{prop:translation.number}
The translation number $\tau : G \to \R$ is independent
of the choice of $x \in \R$ used in its definition.
It is a homomorphism from $G$ to the additive group $\R.$
For any $f \in G$ the set $\Fix(f) \ne \emptyset$
if and only if $\tau(f) = 0.$
\end{prop}

\subsection{Proof of Theorem \ref{theorem:C2interval}}

{\bf Theorem \ref{theorem:C2interval}}
{\em Any nilpotent subgroup of $\Diff^2(\I),\ \Diff^2([0,1))$ or
$\Diff^2(\T^1)$ must be abelian.}

\begin{proof}
Since it suffices to show the action is abelian when restricted to any
invariant interval we may assume our action is irreducible.

\subsubsection{The case $M=\I$ or $M = [0,1)$}

The argument we present here is essentially the same as that given by
Plante and Thurston in (4.5) of \cite{PT}.  We give it here for
completeness and because it is quite simple.

Consider the restriction of $N$ to $(0,1).$  If
no element of $N$ has a fixed point then $N$ is abelian by H\"older's
theorem.  Hence we may assume that there is a non-trivial element $f$
with a fixed point.  Thus, if $h$ is in the center of $N$ Lemma
\ref{lem:fix.component.R} implies that $\Fix(h) \supset \partial \Fix(f)$
which is non-empty.  Another appliction of
Lemma \ref{lem:fix.component.R}, says the non-empty set $\partial \Fix(h)
\subset \Fix(g)$ for any $g \in N.$ We have found a global fixed point
and contradicted the assumption that $N$ acted irreducibly.  We
conclude $N$ is abelian.  \endproof

\subsubsection{The case $M=\T^1$}

Let $N<\Diff^2(\T^1)$ be nilpotent.  As we observed in Theorem
\ref{thm:inv.measure} there is an invariant measure $\mu_0$ for the
action.  If any one element has no periodic points then since it is
$C^2$ it is topologically conjugate to an irrational rotation by
Denjoy's theorem (see e.g. \cite{dM}).  Irrational rotations are
uniquely ergodic so the invariant measure must be conjugate to
Lebesgue measure.  It follows that every element is conjugate to a
rotation so $N$ is abelian.  Hence we may assume that every element of
$N$ has periodic points.

For a homeomorphism of the circle with periodic points, every point
which is not periodic is wandering and hence not in the support of any
invariant measure.  We conclude that the periodic points of any
element contain the support $P$ of the measure $\mu_0$.  It follows that
that $P$ is a subset of the periodic points of 
every element of $N$. 

Since $N$ preserves $\mu_0$ the group $\hat N$ of all lifts to $\R$
of elements of $N$ preserves a measure $\mu$ which is the lift of the
measure $\mu_0$ on $\T^1.$  We can use $\mu$ to define the
translation number homomorphism $\tau_{\mu} \hat N \to \R$
as described above.

We observe that commutators of elements in $\hat N$ must have
translation number $0$.  The covering translations in $\hat N$ are in
its center so if $\hat f$ and $\hat g$ are lifts of $f$ and $g$
respectively then $[\hat f, \hat g]$ is a lift of $[f,g]$.  Since
$\tau_\mu([\hat f, \hat g]) = 0$ implies $[\hat f, \hat g]$ has a
fixed point we may conclude that every element of the commutator
subgroup $N_1 = [N,N]$ must have a fixed point.  Hence any commutator
fixes every element of $P$ because it fixes one point and hence all
its periodic points are fixed.

If $N$ is not abelian there are $f$ and $g$ in $N$ such that $h =
[f,g]$ is a non-trivial element of the center of $N$.  Conjugating
the equation $f^{-1}g^{-1}fg = h$ by $f$ gives $g^{-1}fgf^{-1} = h$,
so $fgf^{-1} = hg$.  Repeatedly conjugating by $f$ gives
$f^ngf^{-n} = h^ng$.  From this we get
$g^{-1}f^ng = h^n f^n$ and by repeatedly conjugating with $g$ we obtain
$g^{-m}f^ng^m = h^{mn} f^n$.  We
conclude $g^{-m}f^ng^mf^{-n} = [g^m, f^{-n}] = h^{mn}.$
Now let $x \in P$ and let $m$ and $n$ be its period under the
maps $f$ and $g$ respectively.  Then $h$,\ $f^n$ and $g^m$ all fix
the point $x$.  If we split the circle at $x$ we get an interval
and a $C^2$ nilpotent group of diffeomorphisms 
generated by $h,\ f^n$ and $g^m$.
By Theorem \ref{theorem:C2interval} this group is abelian. 
Hence $[g^m, f^{-n}] = h^{mn}$ is the identity.  But the only finite
order orientation-preservng homeomorphism of an inteval is the
identity.  We have contrdicted the assumption that $N$ is not abelian.
\end{proof}

\subsection{The proof of  Theorem \ref{thm:fixed.abelian}}

{\bf Theorem \ref{thm:fixed.abelian}}
{\em If $N$ is a nilpotent
subgroup of $\Diff^2(\R)$ and every element of $N$ has a fixed
point then $N$ is abelian.}

\begin{proof}
We first show that there is a global fixed point for $N$.  Consider
$\Fix(h)$ the fixed point set of some non-trivial element $h$ of the
center of $N$. We can apply Lemma \ref{lem:fix.component.R}, and
observe that for any $f \in N, \ \partial \Fix(h) \subset \Fix(f).$
Thus the non-empty set $\partial \Fix(h)$ is fixed pointwise by every
element of $N$.

If $a \in \Fix(h)$ then both $[a,\infty)$ and $(-\infty,a]$ are
invariant by $N$ and each is diffeomorphic to $[0,1).$ We can thus
apply Theorem \ref{theorem:C2interval} and conclude the restriction of
$N$ to each of them is abelian.
\end{proof}

\subsection{The proof of  Theorem \ref{thm:derived.2}}
{\bf Theorem \ref{thm:derived.2}}
{\em Every nilpotent 
subgroup of $\Diff^2(\R)$ is metabelian, i.e.\ has abelian 
commutator subgroup.}

\begin{proof}
Let $N<\Diff^2(\R)$ be nilpotent.  By Theorem \ref{thm:inv.measure} there
is an invariant Borel measure $\mu$ on $\R$ which is finite on 
compact sets and from it we can define a translation number homomorphism
$\tau_\mu : N \to \R.$  According to Proposition \ref{prop:translation.number}
If $N_0$ is the kernel of $\tau_\mu$ then every element of $N_0$ has a fixed
point.  Hence the result follows from Theorem \ref{thm:fixed.abelian}.
\end{proof}

\section{PL actions}
\label{section:PL}

Recall that a {\em piecewise linear homeomorphism} of a one-manifold $M$
is a homeomorphism $f$ for which there exist finitely many subintervals
of $M$ on which $f$ is linear.  Let $\PL(M)$ denote the group of
piecewise-linear homeomorphisms of $M$.

\begin{theorem}
Any nilpotent subgroup of $\PL(\I)$ or $\PL(\T^1)$ is abelian. 
\end{theorem}

\begin{proof}
Suppose first that $N<\PL(\I)$.  Note that there is a natural
homomorphism 
$\psi:\PL(\I)\rightarrow \R^\ast\times \R^\ast$
given by 
$$\psi(f)=(f'(0),f'(1))$$
As $\R\times \R$ is abelian $\psi(f)=(1,1)$ for any nontrivial
commutator in $\PL(\I)$. 

If $N$ is not abelian then there is a nontrivial commutator $f$ lying in
the center of $N$.  Hence $f'(0) = 1$, so that $f$ is the identity on
some interval $[0,a]$ which we take to be maximal. Since every element
of $N$ commutes with $f$, we must have that $g([0,a]) = [0,a]$ for every
$g \in N.$ But then $a$ is also a fixed point of every element of $N$,
so that $f$ is a nontrivial commutator of elements which fix $a$.  We
conclude that $f'(a) = 1$ also.  Hence $f$ is the identity on $[a, b]$
for some $b > a$ and we contradict the maximality of $a$.  We conclude
that the center of $N$ has no nontrivial commutator, so $N$ is
abelian.

Now if $N<\PL(\T^1)$ then the proof is nearly word for word
identical to the proof of Theorem \ref{theorem:C^2_T^1}.  We note
that Denjoy's Theorem (see \cite{dM}) is also valid for PL
homeomorphisms of the circle.  That is, any PL homeomorphism of $\T^1$
without periodic points is topologically conjugate to an irrational
rotation.  With that observation the proof that a PL action of $N$ on $\T^1$ is
abelian is identical to the proof of Theorem \ref{theorem:C^2_T^1}
except that the appeal to Theorem \ref{theorem:C2interval} in the case 
of the interval is
replaced by an appeal to the first part of this theorem, namely
the fact that a PL action of a nilpotent group on $\I$ is abelian.
\end{proof}

\bigskip
\noindent
Benson Farb:\\
Dept. of Mathematics, University of Chicago\\
5734 University Ave.\\
Chicago, Il 60637\\
E-mail: farb@math.uchicago.edu
\medskip

\noindent
John Franks:\\
Dept. of Mathematics, Northwestern University\\
Evanston, IL 60208\\
E-mail: john@math.northwestern.edu


\begin{thebibliography}{ABCDEF}

\bibitem[FF1]{FF1}
B. Farb and J. Franks, Groups of homeomorphisms of one-manifolds, I:
actions of nonlinear groups , preprint, June 2001.

\bibitem[FF2]{FF2}
B. Farb and J. Franks,
Group actions on one-manifolds, II: Extensions of H\"{o}lder's
Theorem, preprint, April 2001.

\bibitem[FS]{FS}
B. Farb and P. Shalen, Groups of real-analytic diffeomorphisms of the
circle, to appear in {\em Ergodic Theory and Dynam. Syst.}

\bibitem[Gh]{Gh}
E. Ghys, Sur les Groupes Engendr\`{e}s par des Diff\`{e}omorphismes
Proches de l'Identit\`{e}, {\em Bol. Soc. Bras. Mat.}, Vol.24, No.2,
p.137-178. 

\bibitem[K]{K}
N. Kopell, Commuting Diffeomorphisms
in {\it Global Analysis, Proceedings of
the Symposium on Pure Mathematics} {\bf XIV}, Amer. Math. Soc.,
Providence, RI (1970), 165--184.

\bibitem[P]{P}
J. Plante, Solvable Groups acting on the line,
Trans. Amer. Math. Soc. {\bf 278}, (1983) 401--414.

\bibitem[PT]{PT}
J. Plante and W. Thurston, Polynomial Growth in Holonomy Groups of Foliations
Comment. Math. Helvetica 39 {\bf 51} 567-584.


\bibitem[dM]{dM}
Welington de Melo,
{\em Lectures on One-dimensional Dynamics} Intstituto de Matematica
Pura e Aplicada do CNPq, Rio de Janeiro, (1989) 
p.137-178. 

\bibitem[MR]{MR}
R. Botto Mura and A. Rhemtulla, {\em Orderable groups}, Marcel
Dekker, 1977.

\bibitem[Ra]{Ra}
M.S. Raghunathan, {\em Discrete subgroups of Lie groups}, Ergebnisse der 
Mathematik, Band 68. Springer-Verlag, 1972. 

\bibitem[Wi]{Wi}
D. Witte, Arithmetic groups of higher $\Q$-rank cannot act on
1-manifolds, {\em Proc. AMS}, Vol.122, No.2, Oct. 1994, pp.333-340.

\end{thebibliography}
\end{document}